\documentclass[final,12p,times,sort&compress]{elsarticle}
\usepackage{graphicx}
\usepackage{epstopdf}
\usepackage{amsmath}
\usepackage{verbatim}
\usepackage{color}
\usepackage{subfigure}
\usepackage{amssymb}

\newcommand{\FT}[1]{\widehat{#1}}

\newcommand{\eproof}{\hfill\rule{2.2mm}{3.0mm}}
\newcommand{\Proof}{\noindent {\bf Proof.~~}}

\newcommand{\R}{{\mathbb R}}
\newcommand{\C}{{\mathbb C}}
\newcommand{\E}{{\mathbb E}\,}
\newcommand{\PP}{{\mathbb P}}
\renewcommand{\eqref}[1]{(\ref{#1})}
\newcommand{\inner}[1]{\langle #1 \rangle}
\newcommand{\diag}{{\rm diag}}
\newcommand{\dist}{{\rm dist}}
\newcommand{\re}{{\rm Re}}
\newcommand{\im}{{\rm Im}}
\newcommand{\vx}{{\mathbf x}}

\newcommand{\vz}{{\mathbf z}}
\newcommand{\va}{{\mathbf a}}
\newcommand{\vv}{{\mathbf v}}
\newcommand{\vu}{{\mathbf u}}
\newcommand{\vf}{{\mathbf f}}
\newcommand{\vg}{{\mathbf g}}
\newcommand{\vs}{{\mathbf s}}
\newcommand{\ve}{{\mathbf e}}
\newcommand{\vh}{{\mathbf h}}
\renewcommand{\H}{{\mathbb F}}
\newcommand{\HH}{{\H}^d}
\newcommand{\NN}{{\mathcal N}}
\newcommand{\F}{{\mathcal F}}
\newcommand{\0}{\mathbf 0}
\renewcommand{\SS}{{\mathcal S}}

\newtheorem{prop}{Proposition}[section]
\newtheorem{lem}[prop]{Lemma}
\newtheorem{defi}{Definition}[section]
\newtheorem{theo}[prop]{Theorem}

\begin{document}

\begin{frontmatter}

\title{Phase retrieval for sub-Gaussian measurements}
\author[label1]{Bing Gao}
\ead{gaobing@naikai.edu.cn}
\author[label2,label3]{Haixia Liu\corref{cor1}}
\ead{liuhaixia@hust.edu.cn}
\author[label4]{Yang Wang}
\ead{yangwang@ust.hk}
\address[label1]{School of Mathematical Sciences, Nankai University, Tianjin 300071, China.}
\address[label2]{School of Mathematics and Statistics, Huazhong University of Science and Technology, Wuhan 430074, China.}
\address[label3]{Hubei Key Laboratory of Engineering Modeling and Scientific Computing, Huazhong University of Science and Technology, Wuhan 430074, China.}
\address[label4]{Department of Mathematics, The Hong Kong University of Science and Technology, 
Clear Water Bay, Kowloon, Hong Kong.}
\cortext[cor1]{Corresponding author.}
\footnote{Authors are alphabetically ordered. All authors have contributed equally to the research.}

\begin{abstract}
Generally, phase retrieval problem can be viewed as the reconstruction of a function/signal from only the magnitude of the linear measurements. These measurements can be, for example, the Fourier transform of the density function. Computationally the phase retrieval problem is very challenging. Many algorithms for phase retrieval are based on i.i.d. Gaussian random measurements. However, Gaussian random measurements remain one of the very few classes of measurements. In this paper, we develop an efficient phase retrieval algorithm for sub-gaussian random frames. We provide a general condition for measurements and develop a modified spectral initialization. In the algorithm, we first obtain a good approximation of the solution through the initialization, and from there we use Wirtinger Flow to solve for the solution. We prove that the algorithm converges to the global minimizer linearly.

\end{abstract}


\end{frontmatter}

\section{Introduction}
\setcounter{equation}{0}

The classic phase retrieval problem concerns the reconstruction of a function from the magnitude of its Fourier transform. Let $f(x) \in L^2(\R^d)$. It is well known that $f$ can be uniquely reconstructed from $\FT f$, where $\FT f$ denotes the Fourier transform of $f$. In many applications such as X-ray crystallography, however, we can only measure the magnitude $|\FT f|$ of the Fourier transform while the phase information is lost. This raises the question whether reconstruction of $f$ (namely recovery of the lost phase information) is possible, up to some obvious ambiguities such as translation and reflection.

Recent focus has been largely on the finite dimensional generalization of the phase retrieval problem. In this setting, one aims to recover a real or complex vector (signal) $\vx$ from the magnitude of some linear measurements of $\vx$. Our paper studies phase retrieval in this setting. On the finite dimensional space  $\H^d$ where $\H=\R$ or $\H=\C$, a set of elements $\F=\{\vf_1,\ldots,\vf_N\}$ in $\H^d$ is called a {\em frame} if it spans $\H^d$. Given this frame, any vector $\vx\in {\H^d}$ can be reconstructed from the inner products $\{\inner{\vx,\vf_j}\}_{j=1}^N$. Often it is convenient to identify the frame $\F$ with the corresponding {\em frame matrix} $F=[\vf_1,\vf_2, \dots, \vf_N]$. The phase retrieval problem in $\HH$ is:

\medskip
\noindent
{\bf The Phase Retrieval Problem.}~~{\em Let $\F=\{\vf_1,\ldots,\vf_N\}$ be a frame in $\H^d$. Can we reconstruct any $\vx\in\HH$ up to a unimodular scalar from $\{|\inner{\vx,\vf_j}|\}_{j=1}^N$, and if so, how?}

\medskip

$\F$ is said to have the {\em phase retrieval (PR) property} if the answer is affirmative. The above phase retrieval problem has important applications in imaging, optics, communication, audio signal processing and more \cite{chai2010array,harrison1993phase,heinosaari2013quantum,millane1990phase,walther1963question}. One of the many challenges is the ``how'' part of the problem, namely to find robust and efficient algorithms for phase retrieval. This turns out to be much more difficult than it looks.   

The phase retrieval problem is an example of a more general problem: the recovery of a vector $\vx\in \H^d$ from {\em quadratic} measurements. For this problem we would like  to recover a vector $\vx\in\H^d$ from a finite number of quadratic measurements $\{\vx^* A_j \vx\}_{j=1}^N$ where each $A_j$ is a Hermitian matrix in $\H^{d\times d}$. This is the so-called {\em generalized phase retrieval problem}, which was first studied in \cite{wang2017generalized} from a theoretical angle, but earlier in special cases such as that for orthogonal projection matrices $\{A_j\}_{j=1}^N$ by others \cite{edidin2017projections,heinosaari2013quantum,cahill2013phase}. 

To computationally recover the signal in phase retrieval, the greatest challenge comes from the nonconvexity of the objective function when it is phrased as an optimization problem.  Let $\F=\{\vf_j\}_{j=1}^N$ in $\H^d$ be the measurement frame for the phase retrieval problem. Assume that  $|\inner{\vx,\vf_j}|^2 = y_j$. A typical set up is to solve the optimization problem 
\begin{equation}  \label{PR-optimization}
        \hat\vx = \mathop{\rm argmin}\limits_{\vx\in\H^d} \frac{1}{N}\sum_{j=1}^N \big(|\inner{\vx,\vf_j}|^2- y_j\big)^2.
\end{equation}
Clearly here the objective function $E(\vx):= \frac{1}{N}\sum_{j=1}^N \big(|\inner{\vx,\vf_j}|^2- y_j\big)^2$ is nonconvex. The same holds for other objective functions used for phase retrieval. As a result, for a general frame, finding the global minimizer of  the optimization problem (\ref{PR-optimization}) is extremely challenging if not intractable.

Nevertheless one class of phase retrieval problems for which very efficient reconstructive algorithms have been extensively studied is when the measurements are i.i.d. Gaussian random measurements. Several approaches based on convex relaxation techniques, such as {\em PhaseLift} \cite{candes2014solving}, {\em PhaseCut} and  {\em MaxCut } have been developed, see \cite{candes2015phase1,waldspurger2015phase}, {\em PhaseMax} \cite{goldstein2018phasemax} and the work by Bahmani and Romberg \cite{bahmani2016phase}. Such convex methods can be computationally challenging for large dimensional problems or high computational complexity, which had led to the development of various non-convex optimization approaches. 
 The methods by {\em AltMinPhase} \cite{netrapalli2013phase} and {\em Karczmarz }\cite{wei2015solving} first estimate the missing phase information and solve the phase retrieval problem through the least square method and Karczmarz method, respectively. It is shown that AltMinPhase converges linearly to the true solution up to a unimodular scalar.
The Wirtinger Flow (WF) algorithm introduced in \cite{candes2015phase} is guaranteed to converge linearly to the global minimizer for Gaussian measurements when the number of measurements $N$ is in the order of $O(d\log d)$. Various other techniques, such as truncated methods \cite{chen2015solving,wang2017solving}, have been developed to improve its efficiency and robustness with $N=O(d)$ Gaussian measurements. Other techniques, such as Gauss-Newton's method \cite{gao2017phaseless}, rank-1 alternating minimization algorithm \cite{cai2017fast} and composite optimization algorithm \cite{duchi2017solving} have all provided theoretical convergence analysis for Gaussian random measurements. Some of the aforementioned methods such as the WF algorithm also work for  Fourier measurements with a very specially designed random mask, namely the Coded Diffraction model \cite{candes2015phase}. However, those are virtually the only models for which provable fast phase retrieval algorithms have been developed.  {\em In a big picture, the lack of phase retrieval models that go beyond Gaussian measurements is extremely conspicuous. }


The main objective of this paper is to fill the above void by analyzing phase retrieval models for sub-gaussian measurements and developing efficient algorithm for such models. More specifically we consider phase retrieval problems where sub-gaussian random measurements are used instead of the traditional Gaussian measurements. It turns out that this change causes significant more challenge in the analysis due to the lack of rotational symmetry. We overcome the challenge through more refined analysis and a slightly weakened result.

Key to any non-convex methods for phase retrieval is the initialization step, from which an approximation of the true solution is obtained. This approximated solution can then be used to serve as the initial guess for iteration steps to converge to the true solution. Especially, we use Wirtinger Flow as an example, which uses the so-called spectral initialization to obtain an initial guess and then refine the result by gradient descent iterations.
When this initial guess is close enough to the true solution, the gradient descent is guaranteed to converge to the true solution. Spectral initialization or other initialization methods work well for the Gaussian model (and for the admissible Coded Diffraction model), but it fails for general sub-gaussian random measurements models.  That's the reason why we require the random variables in the Coded Diffraction model to be {\em admissible}. Here in this paper we develop a more general spectral initialization that is less stringent than before, and thus can be applied to most sub-gaussian random measurements models and efficiently solve corresponding phase retrieval problem computationally.

Our generalized spectral initialization aims to provide an initial approximation for phase retrieval problem with sub-gaussian random measurements. Consider the phase retrieval problem of recovering $\vx\in\H^d$ from quadratic measurements $\{\vx^* \va_j\va^*_j \vx\}_{j=1}^N$, where  $\va_j, j=1,\ldots,N$ are i.i.d. sub-gaussian random vectors. We will require $\va_j$ to be sampled randomly from a given distribution satisfying certain properties. More precisely, our model requires the following {\em conditions for Generalized Spectral Initialization}:

\vspace{1ex}
\noindent
{\bf Conditions for Generalized Spectral Initialization:} 
\begin{itemize}
\item[(I)]~~Let $\va_j, j=1,\ldots,N$ are i.i.d. sub-gaussian random vectors in $\H^d$ and $A_j=\va_j\va_j^*$. Furthermore with probability one $A_j -\mathbf{D}(A_j)$ is not pure imaginary and $\mathbf{D}(A_j) \neq c_j I$, where $\mathbf{D}(A_j)$ denotes the diagonal matrix corresponding to the diagonal part of $A_j$. 
\item[(II)] ~~There exist constants $\tau_j$ independent of $\vx$, $1 \leq j \leq 4$, such that $\E(\va_j)=\0$, $\E(A_j) = \tau_1 I$, and
\begin{align}  \label{gen-SI}
&\E\left((\vx^*A_j\vx) A_j\right)  =\tau_2\|\vx\|^2I+\tau_3\vx\vx^*+\tau_4\diag\left([|x_1|^2,\ldots,|x_d|^2]\right)
\end{align}
for all $\vx\in\H^d$.
\end{itemize}
%
%
 We shall prove that under this model a good approximation to the true solution of the phase retrieval problem can be obtained provided that $N = O(d\log^2d)$ with $A_j$ satisfying conditions (I) and (II). We also develop an efficient algorithm for solving the phase retrieval problem under this model.

The rest of the paper is organized as follows: In Section 2, we give the generalized spectral initialization and prove that the method can with high probability achieve good initial results provided $ N=O(d\log^2d) $.  In Section 3, we prove that when the measurements satisfy the conditions  (I) and (II), then gradient descent iteration can linearly converge to the global minimizer. Finally, we provide the details of the proofs as well as some auxiliary results in Section \ref{proof-theo-spectral-init} and the Appendix, respectively.

\section{Generalized Spectral Initialization}
\setcounter{equation}{0}

Let $\va_j, j=1,\ldots,N$ be i.i.d. sub-gaussian random vectors satisfying conditions (I) and (II) for generalized spectral initialization and set $A_j:=\va_j\va_j^*$.   Now for any $\vx\in\H^d$ we denote $y_j = \vx^* A_j\vx = |\va_j^*\vx|^2$. The goal of phase retrieval is of course to recover $\vx$ up to a unimodular constant from the measurements $\{y_j\}_{j=1}^N$.  The generalized spectral initialization introduced here aims to provide a good first approximation to $\vx$, and we describe how it works. Define
\begin{equation}  \label{eq2.1}
    Y := \frac{1}{N}\sum_{j=1}^{N}y_j A_j = \frac{1}{N}\sum_{j=1}^{N}(\vx^*A_j\vx)\, A_j.
\end{equation}
Note that
$$
    \E(\vx^*A_j\vx) = \vx^*\E(A_j)\vx = \tau_1 \|\vx\|^2.
$$

\vspace{1ex}

\noindent
{\bf Generalized Spectral Initialization:}~~Let $\{\va_j\}_{j=1}^N$ be i.i.d. sub-gaussian random vectors in $\H^d$ satisfying conditions (I) and (II). Set $A_j:=\va_j\va_j^*$ for $j=1,2, \dots, N$. Let $y_j = \vx^* A_j\vx = |\va_j^*\vx|^2$. Denote $ \rho^2=\frac{1}{\tau_1  N}\sum_{j=1}^{N}y_j$.
Set
\begin{equation}  \label{M-def}
		M=Y-\frac{\tau_4}{\tau_3+\tau_4}\mathbf{D}(Y-\tau_2\rho^2 I),
\end{equation}
where $\mathbf{D}(Y-\tau_2\rho^2 I)$ denotes the diagonal matrix consisting only the diagonal part of matrix $Y-\tau_2\rho^2 I$. 
\begin{defi} 
Let $\vz_0 \in \H^{d}$ be the eigenvector corresponding to the largest eigenvalue of $M$ in (\ref{M-def})  normalized to $ \|\vz_0\|_2 = \rho $. We shall call $\vz_0$ the {\em generalized spectral initialization} for the measurements $\{\va_j\}_{j=1}^N$.
\end{defi}

\vspace{1ex}
We shall show that $\tau_3+\tau_4 >0$ and the vector $\vz_0$ provides a good initial approximation to the true solution  $\vx$ if we have enough measurements, much like the classical spectral initialization for Gaussian measurements.

\begin{lem} \label{lem-basic}
    Let $\{\va_j\}_{j=1}^N$ satisfy conditions (I) and (II). Then we have $\tau_2>0$, $\tau_3>0$ and $\tau_3+\tau_4>0$.
\end{lem}
\Proof  Since all $A_j=\va_j\va_j^*$ are identically distributed we will examine conditions (I) and (II) for $A_1$. Write $A_1=[a_{mn}]$. Taking $\vx=\ve_k$ in (\ref{gen-SI}) yields
$$
 \E(a_{kk}^2) = \tau_2 + \tau_3 +\tau_4, ~~ 
 \E(a_{kk}a_{mm}) = \tau_2\mbox{~if $m\neq k$}, 
 ~~\mbox{and}~ \E(a_{kk}a_{mn}) = 0 \mbox{~if $m\neq n$}.
$$
Since for some $k \neq m$ we have $a_{kk} \neq a_{mm}$ by the assumption that ${\mathbf D}(A_1) \neq c_1 I$ we have in this case $\E^2(a_{kk}a_{mm}) <\E(a_{kk}^2)\,\E(a_{mm}^2)$. It follows that $0<\tau_2<\tau_2+\tau_3+\tau_4$. Thus $\tau_3+\tau_4>0$. 

Now taking $\vx = \ve_k+\ve_m$ with $k \neq m$ and looking at the off diagonal elements in (\ref{gen-SI})  we have 
$$
    \E\big(a_{km}(a_{kk}+a_{mm}+a_{km}+a_{mk})\big) =\E\big(a_{mk}(a_{kk}+a_{mm}+a_{km}+a_{mk})\big) = \tau_3.
$$
It is  easy to see that this yields $2\tau_3 = \E\big((a_{km}+a_{mk})^2\big)$. Since $a_{km}+a_{mk}\in\R$ we must have $\tau_3 \geq 0$. But $A_1-{\mathbf D}(A_1)$ is not pure imaginary, so there must exist $k \neq m$ such that $a_{km}+a_{mk} \neq 0$. It follows that $\tau_3>0$.
\eproof

\vspace{1ex}

\begin{theo}  \label{theo-spectral-init}
	Let $\{\va_j\}_{j=1}^N$ be i.i.d. sub-gaussian random vectors in $\H^d$ satisfying conditions (I) and (II) and set $A_j:=\va_j\va_j^*$. For the phase retrieval problem, given the measurements $y_j =\vx^*A_j\vx, j=1,\ldots,N$ let $\vz_0\in \H^d$ be the corresponding generalized spectral initialization. Then for any $\varepsilon >0$, there exist constants $ c_{\varepsilon}, C_{\varepsilon}>0 $ depending  on $\varepsilon$, such that with probability at least $ 1-1/d^3-2\exp(-c_{\varepsilon} N) $ we have 
\begin{equation}  \label{eq:spectral-init}
   \textup{dist}(\vz_0,\vx)\leq \varepsilon\|\vx\|
\end{equation}
provided $ N\geq C_{\varepsilon}d\log^2 d $. 
\end{theo}
\Proof We shall leave the proof of this theorem to Section \ref{proof-theo-spectral-init}.
\eproof

\vspace{1ex}

The above theorem is a key ingredient for solving the sub-gaussian measurements phase retrieval problem. 
\vspace{1ex}

\section{Phase Retrieval with Sub-Gaussian Random Measurements}
\setcounter{equation}{0}

Throughout this section we shall assume that we have random measurments $ \{\va_j\}_{j=1}^N $ satisfying conditions (I) and (II).  The generalized spectral initialization combined with the Wirtinger Flow (WF) method can solve the phase retrieval with sub-gaussian measurements.

As before and throughout the rest of the paper we denote $A_j = \va_j\va_j^*$. Given  $\vx\in\H^d$ (where $\H=\R$ or $\C$) we have measurements $y_j = \vx^* A_j \vx =|\langle \va_j, \vx\rangle|^2$. To recover $\vx$ we solve the following minimization problem:
\begin{equation}  \label{PR-optimization-gen}
        \hat\vz = \mathop{\rm argmin}\limits_{\vz\in\H^d}\frac{1}{2N} \sum_{j=1}^N (\vz^* A_j \vz- y_j)^2.
\end{equation}
The target function $E_\vx(\vz):= \frac{1}{2N}\sum_{j=1}^N (\vz^* A_j \vz- y_j)^2$ is a 4-th order polynomial and is  not convex. 

 
\begin{defi}
Let $\vx\in\C^d$ be the solution of \eqref{PR-optimization-gen} where $y_j=\vx^*A_j\vx$. For any $\vz\in\C^d$ we define $ \theta(\vz) $ as 
$$
    \theta(\vz): = \underset{\theta\in [0,2\pi)}{\arg\min}\,\|\vz-\vx e^{i\theta}\|.
$$
The {\em distance} between $\vx$ and $\vz$ is defined as
$$
          \dist (\vz,\vx)=\|\vz-\vx e^{i\theta(\vz)}\|.
$$
\end{defi}
We also define the $ \varepsilon $-neighborhood of $ \vx $ by
$$
     \SS(\vx, \varepsilon):=\Bigl\{\vz\in\C^d:\dist(\vz,\vx)\le\varepsilon\|\vx\|\Bigr\}.
$$
%
%

To solve the optimization problem \eqref{PR-optimization-gen} where the measurements $\{A_j=\va_j\va_j^*\}_{j=1}^N$ satisfying conditions (I) and (II), we start from an initial guess $\vz_0$ and iterate via
\begin{equation} \label{grad_descent}
   \vz_{k+1}=\vz_k-\xi\cdot\nabla_{\vz}E_\vx(\vz_k):=\vz_k-\xi\cdot\bigg(\frac{1}{N}\sum^N_{j=1}(\vz_k^*A_j\vz_k-y_j)\,A_j\vz_k\bigg),
\end{equation}
with $ \xi >0$ being the stepsize, where as before $E_\vx(\vz):= \frac{1}{2N}\sum_{j=1}^N (\vz^* A_j \vz- y_j)^2$ is the target function.  We shall show that with proper generalized spectral initialization for the initial guess $\vz_0$ such iterations converge to the global minimizer linearly.

\vspace{1ex}

The linear convergence will follow from the two key lemmas below. From the scaling property of the target function $E_\vx(\vz)$, without loss of generality we may assume the true solution $\vx$ to the optimization problem \eqref{PR-optimization-gen} satisfies $ \|\vx\|_2=1$. Throughout the paper, we adopt the notation $(t)_+: = \max(t, 0)$ and $(t)_-: = \max(-t, 0)$, which represent the positive and negative parts of any $t \in\R$ respectively. Throughout the paper, we use $ c $, $ C $ or subscript forms of them to denote constants, whose value may change from instance to instance but {\em depend only on the sub-gaussian norm} of the distribution of the measurements $\{\va_j\}_{j=1}^{N}$.

\begin{lem}[Local Curvature Condition]	\label{lem_regularity_cond}
	Let $ \vx$ be the solution of the optimization problem \eqref{PR-optimization-gen} with $ \|\vx\|_2=1 $. Assume that the measurement vectors $ \{\va_j\}_{j=1}^N $ satisfy conditions (I) and (II). For any sufficiently small $ \delta>0 $ there exist constants $c, c_\delta, C_\delta>0$ where $c_\delta, C_\delta$ depend on $\delta$, such that for $ N\geq C_{\delta}\,d\log^2d$, with probability greater than $ 1-1/d^3-\exp(-cd)-2\exp(-c_{\delta}N) $ we have
\begin{equation}\label{curvature_condition}
	\re\left(\langle\nabla_{\vz} E_\vx(\vz),\,\vz-\vx e^{i\theta(\vz)}\rangle\right)\ge \frac{\beta-\delta}{4}\cdot\dist^2(\vz,\vx)+\frac{1}{10N}\sum_{j=1}^{N}\big|(\vz-\vx e^{i\theta(\vz)})^*A_j(\vz-\vx e^{i\theta(\vz)})\big|^2
\end{equation}
for all $\vx$ and $\vz\in\SS(\vx,\varepsilon_0)$, where
\begin{equation}\label{epsilon-con}
 \varepsilon_0:=\frac{10}{27\alpha} \left(\sqrt{36|\tau_4|^2+\frac{27\alpha\beta}{10}} -6|\tau_4|\right)
\end{equation}
with $\alpha := \tau_2+\tau_3-(\tau_4)_-$ and $\beta:= \tau_3-(\tau_4)_-$.
\end{lem}
\Proof 
Since  $ \va_j, j=1,\ldots,N $ are i.i.d. sub-gaussian random vectors, we may without loss of generality assume $ \max_j \|\va_j\|_{\psi_2}=1 $. By the definition of sub-gaussian random vectors, with probability greater than $ 1-\exp(-cd) $ for some constant $c>0$ we have $ \max_{j}\|\va_j\|\leq \sqrt{ 2d\log N}, j=1,\ldots,N$. 

Let $\vh=e^{-i\theta(\vz)}\vz-\vx$ with $ \|\vh\|\leq \varepsilon_0 $. By definition we have $ \im(\vh^*\vx)=0 $ and 
$$
 \re\left(\langle\nabla_{\vz} E_\vx(\vz),\,\vz-\vx e^{i\theta(\vz)}\rangle\right)
    =\frac{1}{N}\sum^N_{j=1}\left(2\big(\re(\vh^*A_j\vx)\big)^2
    +3\re(\vh^*A_j\vx)(\vh^*A_j\vh)+|\vh^*A_j\vh|^2\right).
$$
To establish (\ref{curvature_condition}), it suffices to prove that
\[
\frac{1}{N}\sum^N_{j=1}\left(2\big(\re(\vh^*A_j\vx)\big)^2
+3\re(\vh^*A_j\vx)(\vh^*A_j\vh)+|\vh^*A_j\vh|^2\right)-\frac{1}{10N}\sum_{j=1}^{N}|\vh^* A_j\vh|^2\ge\frac{\beta-\delta}{4}\|\vh\|^2
\]
holds for all $ \vh $ satisfying $\im(\vh^*\vx)=0 $, $ \|\vh\|\leq \varepsilon_0 $. Equivalently, we only need to prove that for all  $ \vh $ satisfying $\im(\vh^*\vx)=0 $, $ \|\vh\|=1 $ and for all $ s $ with $ 0\le s\le\varepsilon_0 $,
\[
\frac{1}{N}\sum^N_{j=1}\left(2\big(\re(\vh^*A_j\vx)\big)^2
+3s\re(\vh^*A_j\vx)(\vh^*A_j\vh)+\frac{9s^2}{10}|\vh^*A_j\vh|^2\right)\ge\frac{\beta-\delta}{4}.
\]
By Lemma \ref{lem_real}  for $ N\geq C_\delta d\log^2d $, with probability greater than $ 1-1/d^3 - \exp(-c_\delta N) $, we have
\[
\frac{1}{N}\sum^N_{j=1}\big(\re(\vh^*A_j\vx)\big)^2
\le \E\big(\re^2(\vh^*A_j\vx)\big)+\frac{\delta}{2}
\]
for any $ \vh $ with $ \|\vh\|=1 $. Therefore to establish the local curvature condition (\ref{curvature_condition}) it suffices to show that
\begin{equation}\label{curvature_condition1}
	\frac{1}{N}\sum^N_{j=1}\left(\frac{5}{2}\big(\re(\vh^*A_j\vx)\big)^2
	+3s\re(\vh^*A_j\vx)(\vh^*A_j\vh)+\frac{9s^2}{10}|\vh^*A_j\vh|^2\right)\ge\frac{\beta}{4}+\frac{ \E\big(\re^2(\vh^*A_j\vx)\big)}{2}
\end{equation}                                                                                  
To prove this inequality, we first prove it for a fixed $ \vh$, and then use a covering argument. To simplify the statement,  we use the shorthand
\begin{align*}
	U_j : &= \frac{5}{2}\big(\re(\vh^*A_j\vx)\big)^2
	+3s\re(\vh^*A_j\vx)(\vh^*A_j\vh)+\frac{9s^2}{10}|\vh^*A_j\vh|^2\\
	& = \Bigg( \sqrt{\frac{5}{2}}\re(\vh^*A_j\vx)+\sqrt{\frac{9}{10}}s|\vh^*A_j\vh|\Bigg)^2.
\end{align*}      
For a fixed $ \vh $, according to the expectations given in the proof of Lemma \ref{lem_initial_bound} and $ \varepsilon_0<1 $ we have
\begin{align*}
	u_j &= \E U_j\leq \frac{5}{2}\Big(\frac{\tau_3}{2}+\tau_2+\frac{\tau_3}{2}+|\tau_4|\Big)+3\big(\tau_2+\tau_3+|\tau_4|\big)+\frac{9}{10}\big(\tau_2+\tau_3+|\tau_4|\big)\\
	&<7\big(\tau_2+\tau_3+|\tau_4|\big).
\end{align*}                                                                                         
Now define $ X_j=u_j - U_j $. First, since $ U_j\geq 0 $, $ X_j\leq u_j\leq 7\big(\tau_2+\tau_3+|\tau_4|\big) $. Second, we bound $ \E X_j ^2$ using Holder's inequality with $ s\leq \varepsilon_0<1 $:
\begin{align*}
  \E X_j^2\leq \E U_j^2&=\frac{25}{4}\E\big(\re^4(\vh^*A_j\vx)\big)+\frac{81}{100}s^4\E\big(|\vh^*A_j\vh|^4\big)+\frac{27}{2}s^2\E\Big(\re^2(\vh^*A_j\vx)|\vh^*A_j\vh|^2\Big)\\
  &\quad \quad +15s\E\Big(\re^3(\vh^*A_j\vx)|\vh^*A_j\vh|\Big)+\frac{27s^3}{5}\E\Big(\re(\vh^*A_j\vx)|\vh^*A_j\vh|^3\Big)\\
  &\le\frac{25}{4}\sqrt{\E\big(|\va_j^*\vh|^8\big)\E\big(|\va_j^*\vx|^8\big)}+\frac{81}{100}s^4\E\big(|\va_j^*\vh|^8\big)+\frac{27s^2}{2}\sqrt{\E\big(|\va_j^*\vh|^{12}\big)\E\big(|\va_j^*\vx|^4\big)}\\
  &\quad \quad +15s\sqrt{\E\big(|\va_j^*\vh|^{10}\big)\E\big(|\va_j^*\vx|^6\big)}+\frac{27s^3}{5}\sqrt{\E\big(|\va_j^*\vh|^{14}\big)\E\big(|\va_j^*\vx|^2\big)}  \\
  &\leq \bigg(\frac{25}{4}\cdot 8^4+\frac{81}{100}s^4\cdot 8^4 +\frac{27s^2}{2}\cdot \sqrt{12^6\cdot 16}\\
   &\quad \quad+15s\cdot \sqrt{ 6^3\cdot 10^5}+\frac{27s^3}{5}\cdot  \sqrt{ 14^7\cdot 2}\bigg) \cdot \|\va_j\|_{\psi_2}^{8}\\
  &:=C_s.
\end{align*}
Here $ C_s $ is a constant depending only on $ s $ and the second inequality is by the definition of sub-gaussian norm
\begin{align}\label{cond_subgaussian}
 \E(|\va^*\vv|^p) \leq p^{\frac{p}{2}}\|\va\|_{\psi_2}^{p}.
\end{align}
Appling Lemma \ref{one_side_concen} with $ \sigma^2 = N\max\big\{49\big(\tau_2+\tau_3+|\tau_4|\big)^2, C_s \big\} $ and $ y = \frac{\beta}{8}N $,
\[
\PP\bigg(Nu_j - \sum_{j=1}^{N}U_j\geq\frac{\beta}{8}N \bigg)\leq e^{-3\gamma N}.
\]
Therefore, with probability at least $ 1-\exp(-3\gamma N) $, we have
\begin{align*}
	\frac{1}{N}\sum_{j=1}^{N}U_j&\geq u_j- \frac{\beta}{8}\\
	&=\frac{1}{2}\E\big(\re^2(\vh^*A_j\vx)\big)+2\E(\re^2(\vh^*A_j\vx)) + 3s\E(\re(\vh^*A_j\vx)|\vh^*A_j\vh|)+\frac{9}{10}s^2\E(|\vh^*A_j\vh|^2)-\frac{\beta}{8}\\
	&\geq \frac{1}{2}\E\big(\re^2(\vh^*A_j\vx)\big)+\frac{\beta}{2}-\frac{\beta}{8}\\
	&\geq \frac{1}{2}\E\big(\re^2(\vh^*A_j\vx)\big)+\frac{3\beta}{8}.
\end{align*}
Here the second inequality comes from Lemma \ref{lem_initial_bound}.

The inequality above holds for a fixed $ \vh $ and a fixed value $ s $. To prove (\ref{curvature_condition1}) for all $ s\leq\varepsilon_0 $ and all $ \vh\in\C^d $ with $\|\vh\|=1 $, define
\[
p_j(\vh,s)=\sqrt{\frac{5}{2}}\re(\vh^*A_j \vx)+\sqrt{\frac{9}{10}}s|\vh^*A_j\vh|
\]
and
\[
p(\vh,s)=\frac{1}{N}\sum_{j=1}^{N}p_j^2(\vh,s).
\]
Recall that $ \max_{j\in[N]}\|\va_j\|\leq \sqrt{ 2d\log N}, j=1,\ldots,N$ and $ s<1 $, we have $ |p_j(\vh,s)|<6d\log N$. Moreover, for any unit vectors $ \vu, \vv\in\C^d $,
\[
|p_j(\vu,s)-p_j(\vv,s)|\leq\sqrt{\frac{5}{2}}\Big|\re\big((\vu-\vv)^*A_j\vx\big)\Big|+\sqrt{\frac{9}{10}}s\Big(|\va_j^*\vu|+|\va_j^*\vv|\Big)|\va_j^*(\vu-\vv)|<8d\log N \|\vu-\vv\|.
\]
So we have
\begin{align*}
 |p(\vu,s )-p(\vv,s)|&=\Big|\frac{1}{N}\sum_{j=1}^{N}\big(p_j^2(\vu, s)-p_j^2(\vv,s)\big)\Big|\\
 &=\Big|\frac{1}{N}\sum_{j=1}^{N}\big(p_j(\vu, s)+p_j(\vv,s)\big)\big(p_j(\vu, s)-p_j(\vv,s)\big)\Big|\\
 &<96d^2\log^2 N\|\vu-\vv\|.
\end{align*}
Thus when $ \|\vu-\vv\|\leq\eta:=\frac{\beta}{1536d^2\log^2 N} $,
\begin{equation}\label{covering_cond1}
p(\vv,s)\ge p(\vu,s)-\frac{\beta}{16}.
\end{equation}
Let $ \NN_\eta $ be an $ \eta- $net for the unit sphere of $ \C^d $ with cardinality obeying $ |\NN_\eta|\leq(1+2/\eta)^{2d} $. Applying (\ref{curvature_condition1}) together with the union bound, we conclude that for all $ \vu\in\NN_\eta $ and a fixed $ s $,
\begin{align}\label{covering_cond2}
  \PP\,\Bigg(p(\vu,s)&\geq  \frac{1}{2}\E\big(\re^2(\vu^*A_j\vx)\big)+\frac{3\beta}{8}\Bigg)\ge 1-|\NN_\eta|\exp(-3\gamma N)\\\nonumber
  &\geq 1-(1+3072d^2\log^2 N/\beta)^{2d}\exp(-3\gamma N)\\\nonumber
  &\geq 1-\exp(-2\gamma N).
\end{align}
The last line follows by choosing $ N $ as before such that $ N\geq C d\log^2 d $, where $ C $ is a sufficiently large constant. Now for any $ \vh $ on the unit sphere of $ \C^d $, there exists a vector $ \vu\in\NN_\eta $ such that $ \|\vh-\vu\|\leq \eta $. By combining (\ref{covering_cond1}) and (\ref{covering_cond2}), $ p(\vh,s)\geq  \E\big(\re^2(\vh^*A\vx)\big)/2+5\beta/16$ holds with probability at least $ 1-\exp(-2\gamma N) $ for all $ \vh $ and for a fixed $ s $. Applying a similar covering number argument over $ s\leq \varepsilon_0 $ we can further conclude that for all $ \vh $ and $ s $, 
$$ 
p(\vh,s)\geq  \frac{\E\big(\re^2(\vh^*A\vx)\big)}{2}+\frac{\beta}{4} 
$$ 
holds with probability at least $ 1-\exp(-\gamma N) $ as long as $ N\geq C d\log^2 d $. Thus when $ N\geq C_\delta d\log^2d $, (\ref{curvature_condition}) holds with probability greater than $ 1-1/d^3-\exp(-cd)-2\exp(-c_\delta N) $.
\eproof
\vspace{1ex}

\begin{lem}[Local Smoothness Condition]\label{lem_smooth_cond}
     Under the same assumptions as in Lemma \ref{lem_regularity_cond} and with $\varepsilon_0$ being given in \eqref{epsilon-con}, for $ N\geq C_{\delta}d\log^2 d $  and any $\vz\in\SS(\vx,\varepsilon_0)$, with probability at least $ 1-1/d^3 -\exp(-cd)- \exp(-c_{\delta}\,N) $ we have 
  \begin{equation}	\label{smooth_cond}
      \|\nabla_{\vz} E_\vx(\vz)\|^2
	\leq R\Bigg(\frac{\beta-\delta}{8}\dist^2(\vz, \vx)+\frac{1}{10N}\sum_{j=1}^{N}\big|(\vz-\vx e^{i\theta(\vz)})A_j(\vz-\vx e^{i\theta(\vz)})\big|^2\Bigg).
  \end{equation}
Here $ R $ is defined as
\[
R\geq\max\Bigg(96\frac{\hat{\alpha}^2(1+\delta^2)}{\beta-\delta},\,\, 270\hat{\alpha}(1+\delta)+60d\tau_1(1+\delta)\varepsilon_0^2\Bigg)
\]
with $ \hat{\alpha} := \tau_2+\tau_3+|\tau_4| $, $ \beta:=\tau_3-(\tau_4)_- $ and constants $C_{\delta},c_{\delta}>0$ depending on $\delta $.
\end{lem}
\Proof  Set $\vh:=e^{-i\theta(\vz)}\vz-\vx$. For any $\vu\in\C^d$ with $\|\vu\|=1$, let $\vv=e^{-i\theta(\vz)}\vu$. Recall that $A_j=\va_j\va_j^*$, and we calculate
\begin{align*}
    |\nabla_{\vz} E_\vx(\vz)^*\vu|^2
&=\Big|\frac{1}{N}\sum^N_{j=1}\Big(2\re(\vx^*A_j\vh)+\vh^*A_j\vh\Big)\Big(\vv^*A_j\vx+\vv^*A_j\vh\Big)\Big|^2\\
&=\Big|\frac{1}{N}\sum^N_{j=1}\Big(2\re(\vx^*\va_j\va_j^*\vh)+\vh^*\va_j\va_j^*\vh\Big)\left(\vv^*\va_j\va_j^*\vx+\vv^*\va_j\va_j^*\vh\right)\Big|^2\\
&\le\bigg(\frac{1}{N}\sum^N_{j=1}3|\vx^*\va_j|\,|\vv^*\va_j|\,|\vh^*\va_j|^2+2|\vh^*\va_j|\,|\vv^*\va_j|\,|\vx^*\va_j|^2+|\vh^*\va_j|^3|\vv^*\va_j|\bigg)^2\\
&\le 27\,\Big(\frac{1}{N}\sum^N_{j=1}|\vx^*\va_j|\,|\vv^*\va_j|\,|\vh^*\va_j|^2\Big)^2+12\,\Big(\frac{1}{N}\sum^N_{j=1}|\vh^*\va_j|\,|\vv^*\va_j|\,|\vx^*\va_j|^2\Big)^2+3\,\Big(\frac{1}{N}\sum^N_{j=1}|\vh^*\va_j|^3|\vv^*\va_j|\Big)^2\\
&:= 27I_1+12I_2+3I_3.
\end{align*}
Without loss of generality, we assume that $ \max_j\|\va_j\|_{\psi_2}=1 $. As before the inequality $ \max_{j}\|\va_j\|\leq \sqrt{2d\log N}, j=1,\ldots,N$  holds with probability at least $ 1-\exp(-cd) $. Combined this fact with the
Cauchy-Schwarz inequality and Lemma \ref{concentration_YM} we obtain
\begin{eqnarray*}
	I_1&\le&\bigg(\frac{1}{N}\sum^N_{j=1}|\vh^*\va_j|^4\bigg)\cdot\bigg(\frac{1}{N}\sum^N_{j=1}|\vv^*\va_j|^2|\vx^*\va_j|^2\bigg)\\
	&\le&\frac{1}{N}\sum^N_{j=1}|\vh^*A_j\vh|^2\cdot\vv^*\bigg(\frac{1}{N}\sum^N_{j=1}(\vx^*A_j\vx) A_j\bigg)\vv\\
	&\le&\hat{\alpha}(1+\delta)\cdot\frac{1}{N}\sum^N_{j=1}|\vh^*A_j\vh|^2,
\end{eqnarray*}
\begin{eqnarray*}
	I_2\le\vv^*\bigg(\frac{1}{N}\sum^N_{j=1}(\vx^*A_j\vx) A_j\bigg)\vv\cdot\vh^*\bigg(\frac{1}{N}\sum^N_{j=1}(\vx^*A_j\vx) A_j\bigg)\,\vh\le\hat{\alpha}^2(1+\delta)^2\|\vh\|^2.
\end{eqnarray*}
Moreover
\begin{align*}
	I_3&
	\leq\bigg(\frac{1}{N}\sum^N_{j=1}|\vh^*\va_j|^4\bigg)\cdot\bigg(\frac{1}{N}\sum^N_{j=1}|\vh^*\va_j|^2|\vv^*\va_j|^2\bigg)\\
	&\le\bigg(\frac{1}{N}\sum^N_{j=1}|\vh^*A_j\vh|^2\bigg)\cdot\bigg(\max_{j\in [N]}\|\va_j\|^2\frac{1}{N}\sum^N_{j=1}|\vh^*A_j\vh|\bigg)\\
	&\le \frac{\tau_1(1+\delta)2d\log N}{N}\sum^N_{j=1}|\vh^*A_j\vh|^2\cdot\|\vh\|^2
\end{align*}	
holds with probability at least $ 1-1/d^3-\exp(-cd)-\exp(-c_{\delta}\,N) $ provided  $ N\geq C_{\delta}\,d\log^2d $ for measurements $ A_j, j=1,\ldots, N $ satisfying conditions (I) and (II). 
Therefore, with high probability we have 
\begin{eqnarray*}
	\|\nabla_{\vz} E_\vx(\vz)\|^2&=&\max_{\|\vu\|=1}|\nabla_{\vz} E_\vx(\vz)^*\vu|^2\\
	&\le& 27 \hat{\alpha}(1+\delta)\cdot\frac{1}{N}\sum^N_{j=1}|\vh^*A_j\vh|^2+12\hat{\alpha}^2(1+\delta)^2\|\vh\|^2+6\frac{\tau_1(1+\delta)d\log N}{N}\sum^N_{j=1}|\vh^*A_j\vh|^2\cdot\|\vh\|^2\\
	&=&12\hat{\alpha}^2(1+\delta)^2\|\vh\|^2+\frac{27\hat{\alpha}(1+\delta)+6d\log N\tau_1(1+\delta)\|\vh\|^2}{N}\sum^N_{j=1}|\vh^*A_j\vh|^2\\
	&\le& R\Bigg(\frac{\beta-\delta}{8}\|\vh\|^2+\frac{1}{10N}\sum_{j=1}^{N}|\vh^*A_j \vh|^2\Bigg)
\end{eqnarray*}    
The last line holds as long as
\begin{equation}\label{R_cond}
R=\max\Bigg(96\frac{\hat{\alpha}^2(1+\delta^2)}{\beta-\delta},\,\, 270\hat{\alpha}(1+\delta)+60d\log N\tau_1(1+\delta)\varepsilon_0^2\Bigg).
\end{equation}
\eproof

\begin{prop} \label{conver_rate}
  Under the same assumptions as in Lemma \ref{lem_regularity_cond}, let $\vz_0\in\SS(\vx,\varepsilon_0)$ where  $ \varepsilon_0 $ is given by (\ref{epsilon-con}). Assume that $ N\geq C_\delta d\log^2d $.  Then for $\delta<\beta$ and  $0<\xi\le 2/R$  each iteration $\vz_{k+1}$ given by \eqref{grad_descent} satisfies 
$$
    \dist^2(\vz_{k+1},\vx)\leq\bigg(1-\frac{\xi(\beta-\delta)}{4}\bigg)\cdot\dist^2(\vz_k,\vx)
$$
with probability greater than $ 1-1/d^3 - \exp(-cd)- 2\exp(-c_\delta \,N) $. Here $ \beta = \tau_3-(\tau_4)_- $ and $ R $ defined by (\ref{R_cond}).
\end{prop}
\Proof Assume that $\vz_k \in \SS(\vx,\varepsilon)$ for any $ \varepsilon<\varepsilon_0 $, according to Lemma \ref{lem_regularity_cond} and Lemma \ref{lem_smooth_cond}, for sufficiently small $ \beta>\delta>0 $ when $ N\geq C_\delta d\log^2d $  with probability greater than  $ 1-1/d^3  - \exp(-cd)-2\exp(-c_\delta \,N) $, we have
\[
\re\big(\langle\nabla_{\vz} E_\vx(\vz_k),\,\vz_k-\vx e^{i\theta(\vz_k)}\rangle\big)\geq\frac{\beta-\delta}{8}\|\vz_k-\vx e^{i\theta(\vz_k)}\|^2+\frac{1}{R}\|\nabla_{\vz} E_\vx(\vz_k)\|^2.
\]
Then after one step iteration we obtain
\begin{align*}
\|	&\vz_{k+1}-\vx e^{i\theta(\vz_{k+1})}\|^2\\
	&\le\|\vz_{k+1}-\vx e^{i\theta(\vz_k)}\|^2\\
	&\le\|\vz_k-\vx e^{i\theta(\vz_k)}-\xi\cdot\nabla_{\vz} E_\vx(\vz_k)\|^2\\
	&\le\|\vz_k-\vx e^{i\theta(\vz_k)}\|^2-2\xi\cdot\re\big(\langle\nabla_{\vz} E_\vx(\vz_k),\,\vz_k-\vx e^{i\theta(\vz_k)}\rangle\big)+\xi^2\|\nabla_{\vz} E_\vx(\vz_k)\|^2\\
	&\le\|\vz_k-\vx e^{i\theta(\vz_k)}\|^2-2\xi\cdot\Bigg(\frac{\beta-\delta}{8}\|\vz_k-\vx e^{i\theta(\vz_k)}\|^2+\frac{1}{R}\|\nabla_{\vz} E_\vx(\vz_k)\|^2\Bigg)+\xi^2\|\nabla_{\vz}E_\vx(\vz_k)\|^2\\
	&=\left(1-\frac{\xi(\beta-\delta)}{4}\right)\|\vz_k-\vx e^{i\theta(\vz_k)}\|^2+\xi(\xi-2/R)\|\nabla_{\vz}E_\vx(\vz_k)\|^2\\
	&\leq\left(1-\frac{\xi(\beta-\delta)}{4}\right)\|\vz_k-\vx e^{i\theta(\vz_k)}\|^2.
\end{align*}
The last inequality is according to the fact that $ 0<\xi\le 2/R $.
\eproof

\section{Numerical Implements} \label{numerical}
\setcounter{equation}{0}
We present some numerical experiments to evaluate our proposed algorithm here. In Subsection \ref{numerical_results}, we test the generalized spectral initialization (GSI) and give a comparison with spectral initialization (SI) introduced in \cite{candes2015phase}, and then explore the minimal measurement number we need for successful recovery.

\subsection{Experimental setup}
In the following numerical evaluations we test a couple of different sub-gaussian random measurements. The target signal  $\vx=[x_1,\ldots,x_d]$, which we try to recover, is sampled by 
$$
    x_j=\left\{\begin{array}{cc} \tilde{x}_j,&j<d-1\\
                   200\tilde{x}_j,&j\ge d-1 \end{array} \right.
$$
where $\tilde{\vx}=[\tilde{x}_1,\cdots,\tilde{x}_d]^T\in\H^d$ is random Gaussian. We shall evaluate our method by the relative error of the reconstruction, which is defined as $\dist(\vz,\vx)/\|\vx\|$, 
with $\vz$ being numerical solution and $\vx$ being true solution. Two types of  sub-gaussian random measurements will be used here: uniform measurements and ternary measurements. More precisely, the uniform distribution is  $ \mathbf{U}[-1,1] $ in real case and $ \frac{1}{\sqrt{2}}\mathbf{U}[-1,1]+i\frac{1}{\sqrt{2}}\mathbf{U}[-1,1] $ in complex case.  The ternary distribution we use has the distribution
\begin{equation}\label{ternary_dis}
t= \left\{
\begin{aligned}
1 &\quad \text{with prob.} & 1/3 \\
0 &\quad   \text{with prob.} & 1/3\\
-1 &\quad   \text{with prob.} & 1/3
\end{aligned}
\right.
\end{equation}
in real field and $ \frac{1}{\sqrt{2}}t+i\frac{1}{\sqrt{2}}t $ in complex field. 

For the iteration process we use Barzilai-Borwein's method (B-B method) to choose step size. That is, to use the information in the previous iteration to determine the step size of the current iteration. More precisely, in order to obtain the $(k+1)$-th solution $\vz_{k+1}=\vz_k-\xi_k\nabla_{\vz} E_\vx(\vz_k)$, we may choose
\begin{equation*}
\xi_k=\underset{\xi}{\arg\min} \,\,\|\vs_k-\xi\cdot\vg_k\|^2,
\end{equation*}
where $\vs_k=\vz_k-\vz_{k-1}$ and $\vg_k=\nabla_{\vz} E_\vx(\vz_k)-\nabla_{\vz} E_\vx(\vz_{k-1})$.
By simple calculation, we get
$$
\xi_k=\frac{\re(\langle \vg_k,\,\vs_k\rangle)}{\|\vg_k\|^2}=\frac{ \re(\vg_k^*\vs_k)}{\|\vg_k\|^2}.
$$
In order to ensure step size is positive, we choose the absolute value of $\xi_k$ as the step size. The iterative algorithm stops if $\|\vg_k\|_2<10^{-16}$ or the maximal iteration number is reached.

\subsection{Numerical results}
\label{numerical_results}
\begin{figure}[htbp]
	\begin{center}
		\subfigure[Uniform (real)]{
			\includegraphics[width=0.48\textwidth]{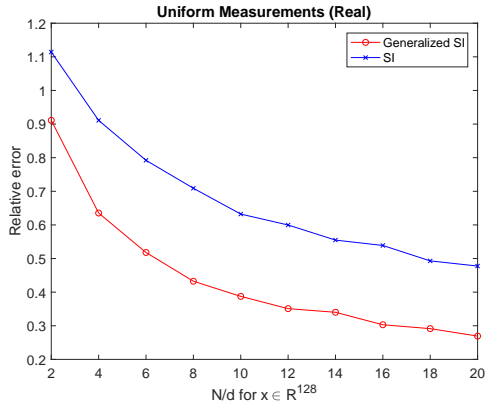}}
		\subfigure[Uniform (complex)]{
			\includegraphics[width=0.48\textwidth]{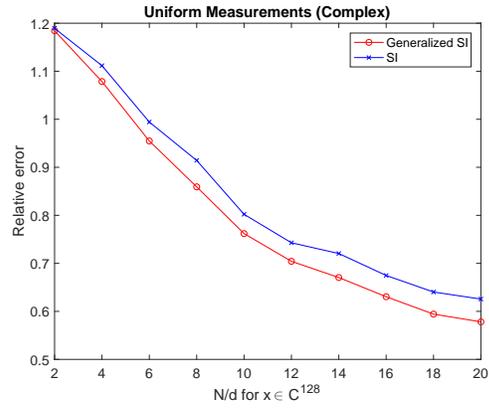}}
		\subfigure[Ternary (real)]{
			\includegraphics[width=0.48\textwidth]{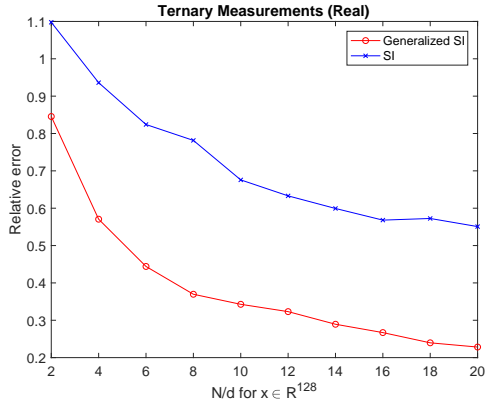}}
		\subfigure[Ternary (complex)]{
			\includegraphics[width=0.48\textwidth]{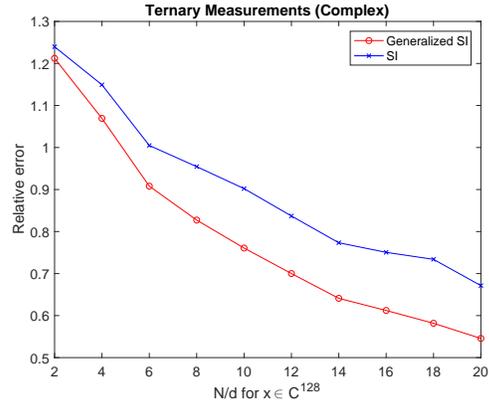}}
		\caption{Initialization experiments: Averaged relative error between $ \vz $ and $ \vx $ for $ d = 128 $ with (a) uniform measurements in the real number field and (b) uniform measurements in the complex number field and (c) ternary measurements in the real number field and (d) ternary measurements in the complex number field} \label{fig:init1}
	\end{center}
\end{figure}

The initialization step of our proposed algorithm is the generalized spectral initialization, which is obtained after 50 iterations of the power method.
First of all, we test GSI by exploring the relationship between the relative error and $N/d$ and give  a comparison with SI. We choose $d=128$ and change $N/d$ in the range $[2,20]$ with stepsize $2$. For each $N/d$, we repeat $50$ times and report the average value of the relative error. Figure \ref{fig:init1} gives the plot of the relative error versus $N/d$ for uniform measurements and  ternary measurements. From Figure \ref{fig:init1}, we observe that under the same measurements, the relative error is much smaller by GSI than that by SI. This is of course expected since SI is designed with only the Gaussian measurements in mind.  While here we only show results for $d=128$, the same conclusion holds for dimensions $d=256,512,1024$ as well in our experiments. 

\begin{figure}[h]
	\begin{center}
		\subfigure[]{
			\includegraphics[width=0.45\textwidth]{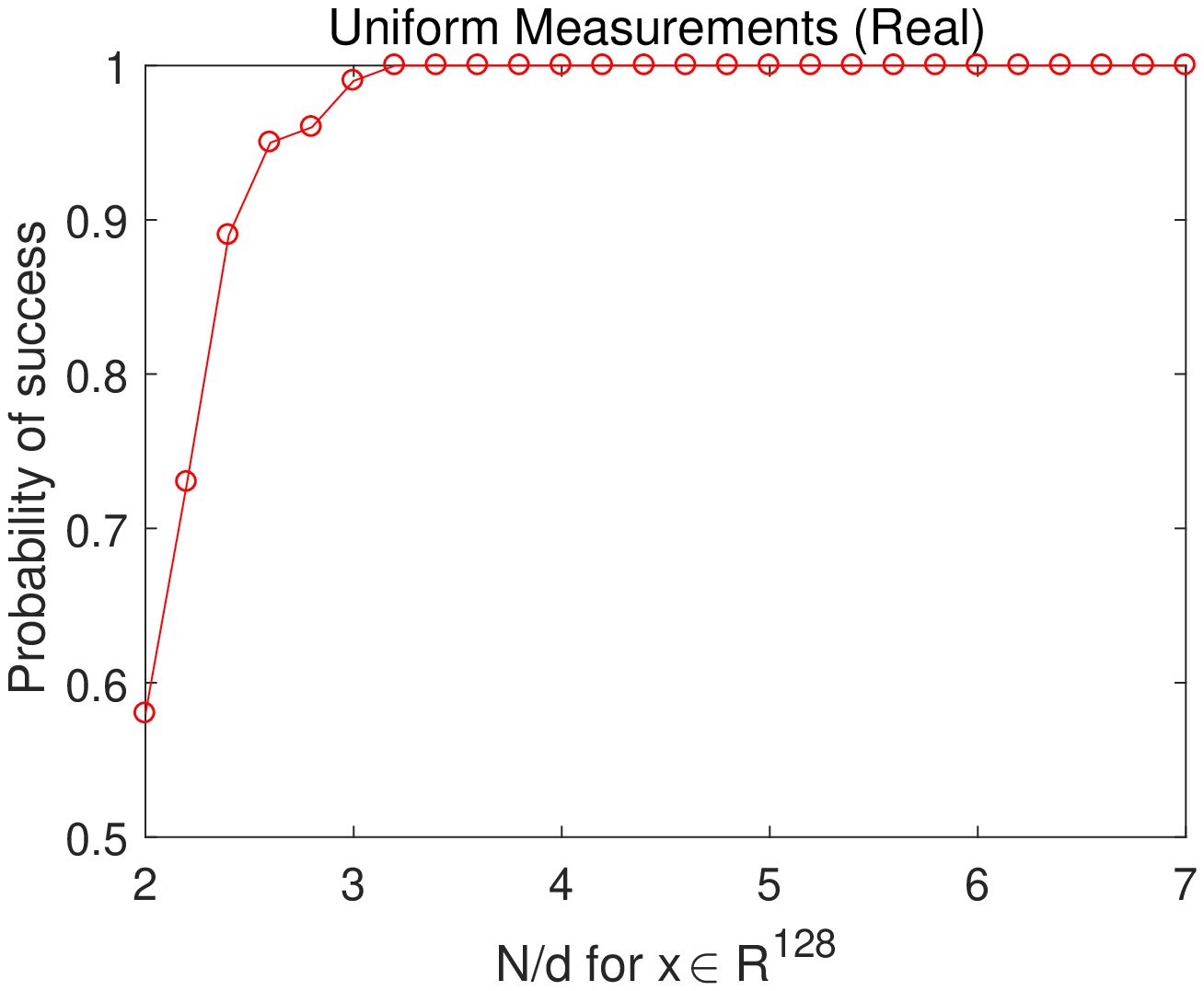}}
		\subfigure[]{
			\includegraphics[width=0.45\textwidth]{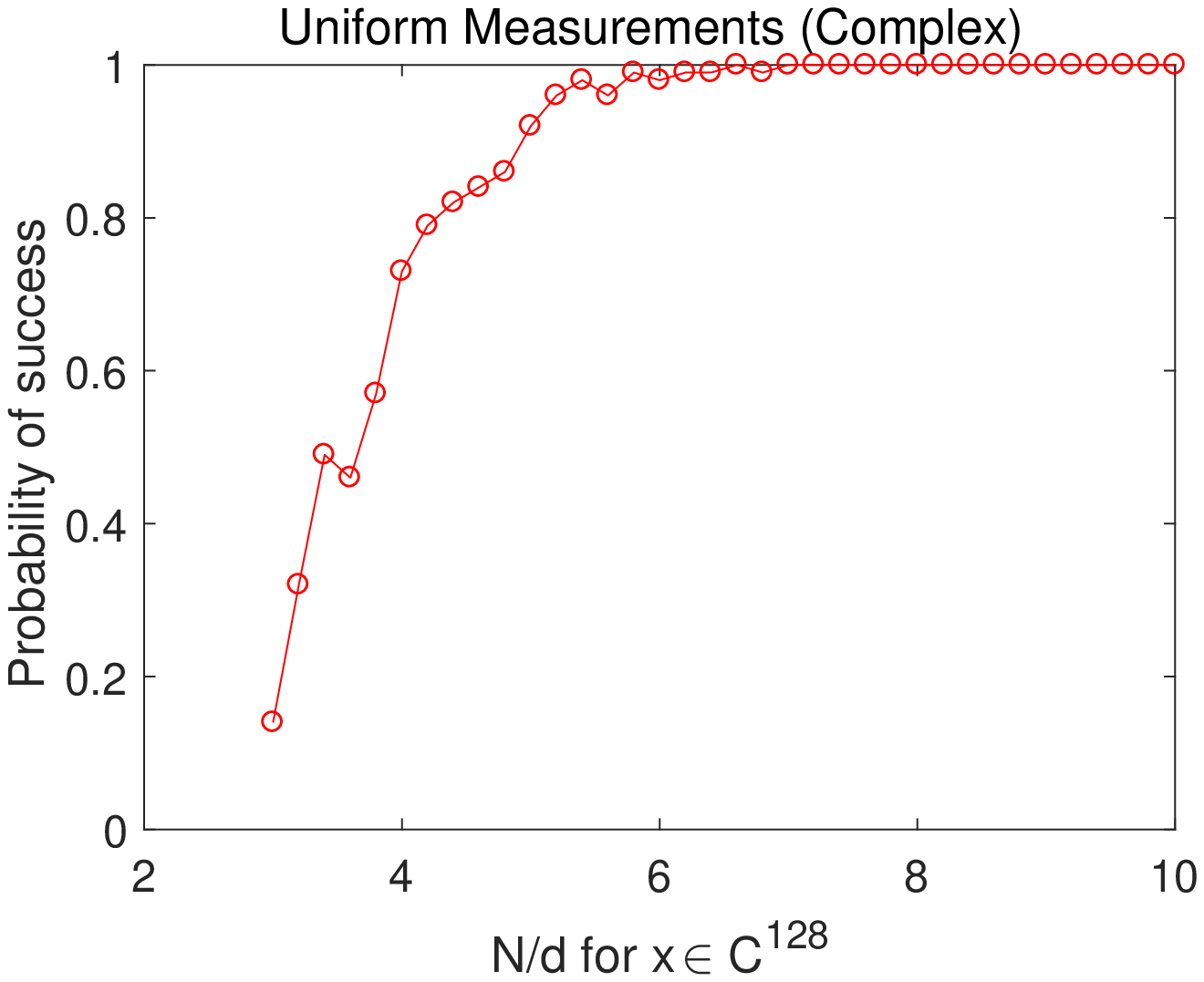}}
		\subfigure[]{
			\includegraphics[width=0.45\textwidth]{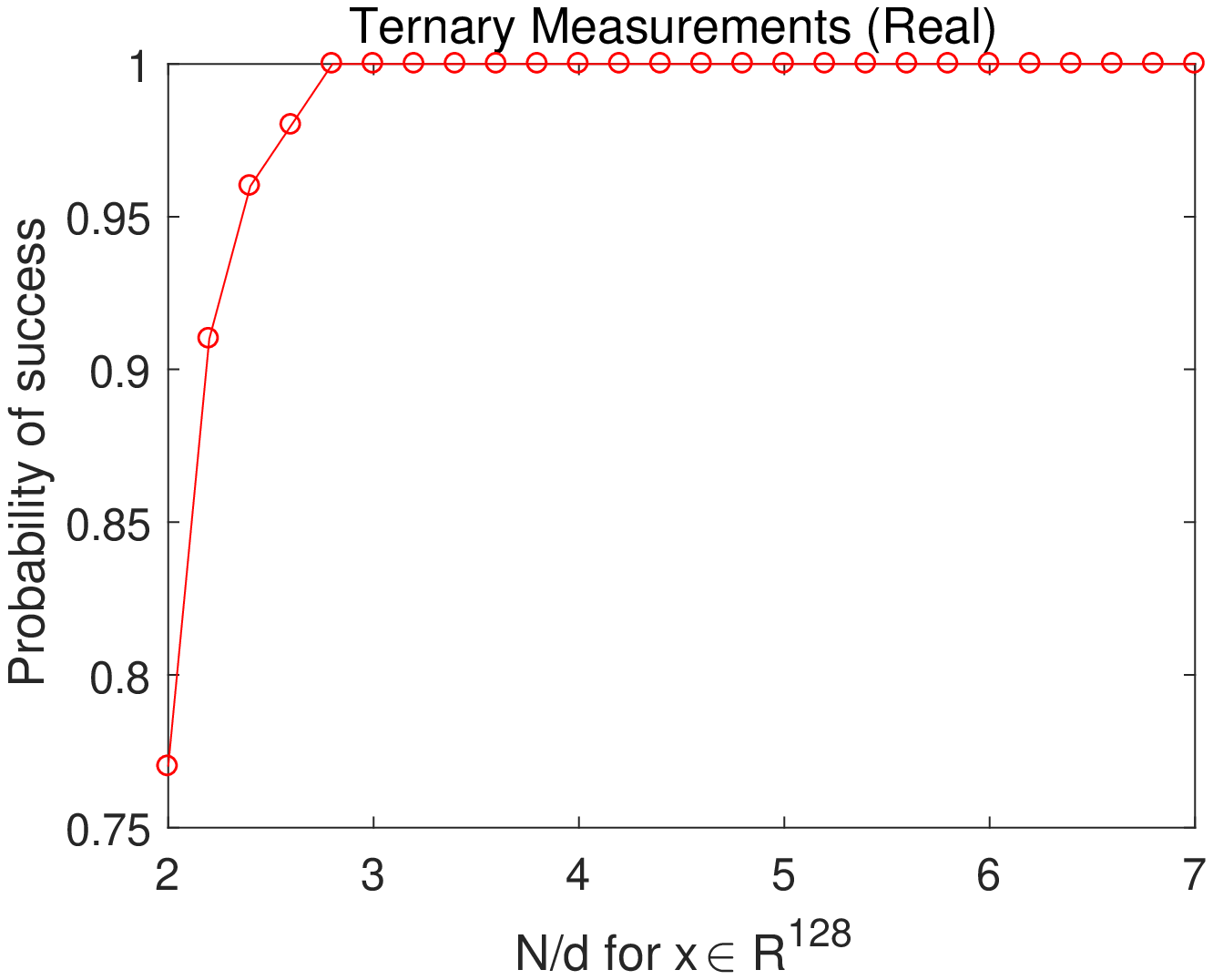}}
		\subfigure[]{
			\includegraphics[width=0.45\textwidth]{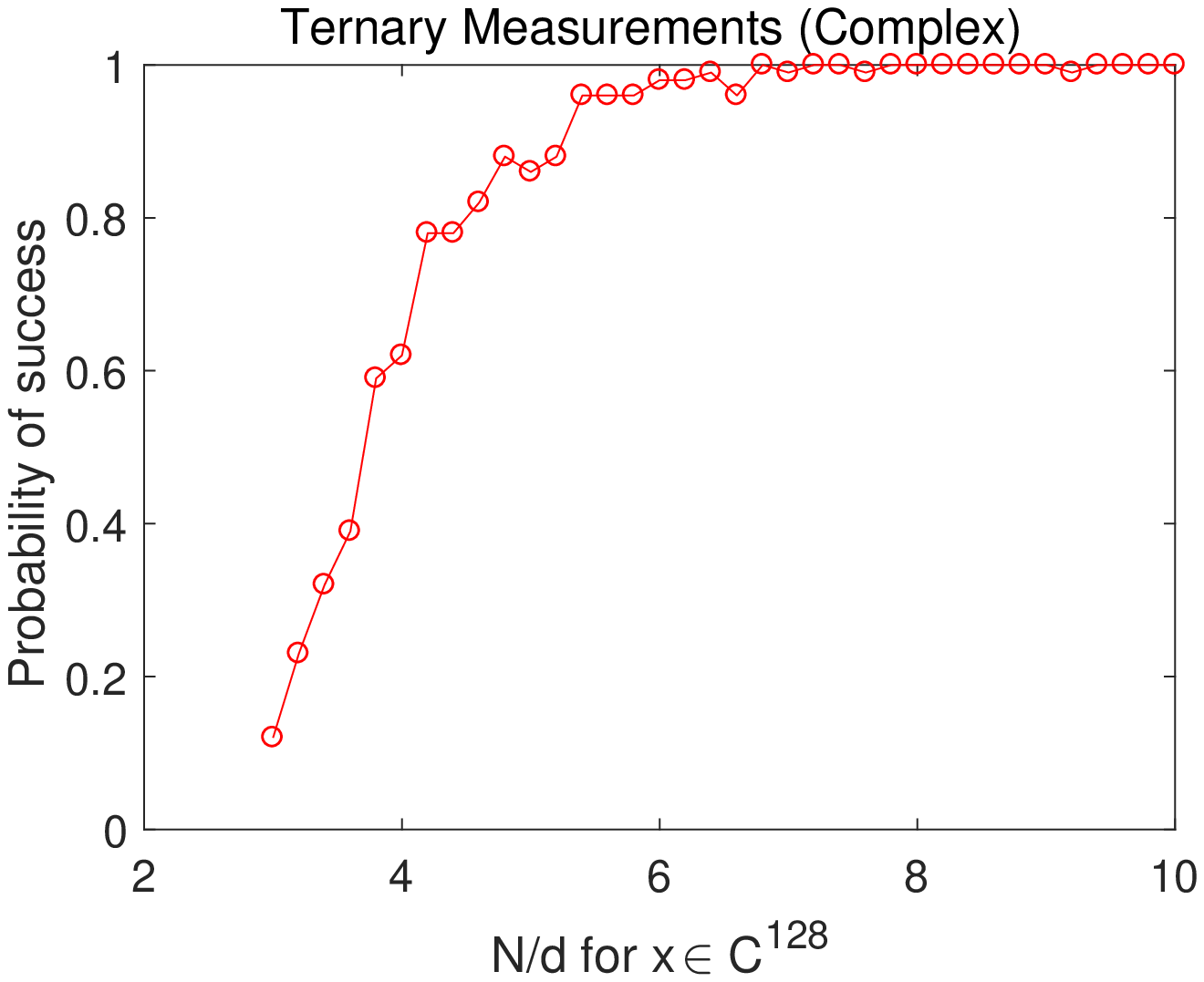}}
		\caption{The plot of success rate versus $N/d$ for uniform and ternary measurements in real and complex field.} \label{fig:prob_success}
	\end{center}
\end{figure}

Starting from the initial guess by generalized spectral initialization, we set the maximal iteration number to 2000 and evaluate our algorithm by 100 trials. In each trial, we declare it a success if the relative error of the reconstruction is less than $10^{-5}$. The empirical probability of success is defined as the average of success rate over 100 trials. 
In our experiments, with $d=128$ Figures \ref{fig:prob_success} shows that $6d$ ($8d$)  measurements in complex field or $3d\ (3d)$ measurements in real field is enough for exact recovery of both the uniform and the ternary measurements with high probability.
%


\section{Proof of Theorem ~\ref{theo-spectral-init}}
\setcounter{equation}{0}
\label{proof-theo-spectral-init}

In this section we prove Theorem \ref{theo-spectral-init} by establishing a series of lemmas.

\begin{lem}[\cite{vershynin2010introduction}, Remark 5.40]
	\label{rmk_2}
	Assume that $ A $ is an $ N\times d  $ matrix whose rows $ A_j$ are independent sub-gaussian random vectors with second moment matrix $ \Sigma $. Let $K$ be the maximum of their sub-gaussian norms. Then for every $t\ge 0$, the following inequality holds with probability at least $1-2\exp(-ct^2/K^4)$:
	$$\left\|\frac{A^*A}{N}-\Sigma\right\|\le\max(\delta,\delta^2),$$
	where $\delta=CK^2\Big(\sqrt{\frac{d}{N}}+\frac{t}{\sqrt{N}}\Big)$ and $C$, $ c $ are absolute constants.
\end{lem}

\begin{lem}\label{norm_x}
	Let $ \{\va_j\}_{j=1}^N $ be sub-gaussian random vectors satisfying conditions (I) and (II) for generalized spectral initialization and $A_j:=\va_j\va_j^* $. Set $ \rho^2=\frac{1}{\tau_1  N}\sum_{j=1}^{N}y_j  $. For any $ \delta>0 $ there exist constants $ C'_\delta, c'_\delta>0$ such that for $ N\geq C'_\delta \,d $, with probability at least $ 1-2\exp(-c'_\delta N) $, we have
	\begin{equation}\label{condition_r}
		|\rho^2-\|\vx\|^2|\leq \delta\|\vx\|^2.
	\end{equation}
\end{lem}
\Proof 	Define
$$
S  = \frac{1}{N}\sum_{j=1}^{N}\vx^* A_j \vx= \frac{1}{N}\sum_{j=1}^{N}\vx^*\va_{j}\va_{j}^*\vx.
$$
Since $ A_j, \,j=1,\ldots,N $ satisfy conditions (I) and (II), we obtain
\[
\E(S) = \tau_1\|\vx\|^2I.
\]
As $ \{\va_j,\,j=1,\ldots,N\} $ are independent  sub-gaussian random vectors,  by Lemma \ref{rmk_2}, for any $ \delta>0 $ there exist constants $ C'_\delta, c'_\delta>0$ such that for $ N\geq C'_\delta \,d $, with probability at least $ 1-2\exp(-c'_\delta N) $, we have
\begin{align*}
|\rho^2-\|\vx\|^2|&= \frac{1}{\tau_1 }|S-\E(S)|\\
&\leq \frac{1}{\tau_1 }\Big|\frac{1}{N}\sum_{j=1}^{N}\vx^*\va_{j}\va_{j}^*\vx-\tau_1\|\vx\|^2\Big|\\
&\leq\frac{1}{\tau_1}\cdot \tau_1\delta\,\|\vx\|^2=\delta\,\|\vx\|^2.
\end{align*}	 
\eproof

\begin{lem}\label{concentration_YM}
    Let $ \{\va_j\}_{j=1}^N $ be sub-gaussian random vectors satisfying conditions (I) and (II) for generalized spectral initialization and $A_j:=\va_j\va_j^* $. Then for  any  $ \delta > 0 $ there exists a $C_{\delta}>0$ such that for $ N\geq C_{\delta} d\log^2 d $, with probability at least $ 1-1/d^3-2\exp(- c_{\delta }N) $, we have
\[
    \|Y-\E(Y)\|\leq \delta\quad \text{and}\quad \|M-\E(M)\|\leq \delta.
\]
\end{lem}
\Proof 
Recall that 
\[
Y = \frac{1}{N}\sum_{j=1}^{N}|\vx^*A_j\vx|\cdot\va_{j}\va_{j}^*.
\]
Define 
\[
\tilde{Y } = \frac{1}{N}\sum_{j=1}^{N}|\vx^*A_j\vx|\cdot \va_{j}\va_{j}^*\cdot I_{\{|\vx^*A_j\vx|\leq2Q\log d\}}.
\]
Here $ Q $ is a positive scalar whose value will be determined shortly. Then define the events: 
\begin{align*}
 E_1(Q)  &= \{\|\tilde{Y}-\E(Y)\|\leq\varepsilon\}\\
 E_2(Q)  &= \{ \tilde{Y} =  Y\}\\
 E_3(Q)  &=\{ |\va_{j}^*\vx|\leq\sqrt{2Q\log d}, j=1,\ldots, N\}\\
 E  &= \{\|Y-\E(Y)\|\leq\varepsilon\}.
\end{align*}
Note that $ E_1(Q) \cap E_2(Q)\subset E $ and $ E_3(Q)\subset E_2(Q) $. Then we have 
\begin{align*}
\PP(E^c)&\leq \PP(E_1^c\cup E_2^c)=\PP\big(E_1^c\cap(E_2\cup E_2^c)\cup E_2^c\big)\\
&=\PP\big((E_1^c\cap E_2)\cup E_2^c\big)\leq \PP(E_1^c\cap E_2)+\PP( E_2^c)\leq \PP(E_1^c\cap E_2)+\PP( E_3^c)
\end{align*}
As for any $ j $, $ \va_{j} $ is a sub-gaussian random vector, we  assume $ \max_j\|\va_j\|_{\psi_2}=1$ without loss of generality. So according to the definition, we obtain
\[
\PP\Big(|\va_j^*\vx|>\sqrt{2Q\log d}\Big)\leq 2d^{-Q}.
\]
This also implies $ \PP(E_3^c)\leq 2N\cdot d^{-Q}$.  Next, we commit to proving that the event $ E_1^c(Q)\cap E_2= \{\|\tilde{Y}-\E(\tilde{Y})\|>\varepsilon\} $ holds with small probability.
Note that $ \vg_j := \sqrt{|\vx^*A_j\vx|}\va_{j}\cdot I_{\{\sqrt{|\vx^*A_j\vx|}\leq\sqrt{2Q\log d}\}}$. Hence $\hat\vg_j:=\frac{ \vg_j}{\sqrt{\log d}} $, $ j=1,\ldots,N $  are i.i.d. sub-gaussian random variables with $ K:=\max_j\|\hat\vg_j\|_{\psi_2} \leq 2Q\|\va_j\|_{\psi_2}$.  Applying Lemma \ref{rmk_2},
\begin{align*}
\Big\|\frac{1}{N}\sum_{j=1}^{N}\hat\vg_j\hat\vg_j^*-\E\Big(\frac{1}{N}\sum_{j=1}^{N}\hat\vg_j\hat\vg_j^*\Big)\Big\|> \frac{\varepsilon}{\log d}
\end{align*}
holds with probability less than $ 2\exp(-c_{\varepsilon,Q}N) $ provided $ N\geq C_{\varepsilon,Q} d\log^2 d $.  This inequality implies 
\begin{align*}
\|\tilde{Y}-\E\tilde{Y}\| &=\Big\|\frac{1}{N}\sum_{j=1}^{N}\vg_j\vg_j^*-\E\Big(\frac{1}{N}\sum_{j=1}^{N}\vg_j\vg_j^*\Big)\Big\|\\
&=\log d\cdot  \Big\|\frac{1}{N}\sum_{j=1}^{N}\hat\vg_j\hat\vg_j^*-\E\Big(\frac{1}{N}\sum_{j=1}^{N}\hat\vg_j\hat\vg_j^*\Big)\Big\|\\
&> \varepsilon.
\end{align*}
Then combined the results, we know
\[
\PP(E^c)\leq 2\exp(-c_{\varepsilon,Q}N)+2N\cdot d^{-Q}
\]
provided $ N\geq C_{\varepsilon,Q} d\log^2 d $. Thus by choosing  $ Q=5 $ and $ N\geq C_{\varepsilon}d \log^2 d $, the bound $ \|Y-\E Y\|\leq \varepsilon$ holds with probability at least $ 1-1/d^3-2\exp(-c_{\varepsilon}N) $.
Note that
\begin{align*}
\|M-E(M)\|&=\Big\|Y-\E(Y)+\frac{\tau_4}{\tau_3+\tau_4}\Big(\E(\mathbf{D}(Y-\tau_2\rho^2 I))-\mathbf{D}(Y-\tau_2\rho^2 I)\Big)\Big\|\\
&\leq \|Y-\E(Y)\|+\frac{|\tau_4|}{\tau_3+\tau_4}\|\E(Y-\tau_2\rho^2 I )-Y+\tau_2\rho^2 I\|\\
&\leq \Big(1+\frac{|\tau_4|}{\tau_3+\tau_4}\Big) \|Y-\E(Y)\|+\frac{\tau_2|\tau_4|}{\tau_3+\tau_4}\big|\rho^2-\E(\rho^2)\big|
\end{align*}
Then combining the above conclusion and Lemma \ref{norm_x}, for any $ \delta>0 $, we choose $ \varepsilon=\frac{\delta}{2\tau} $ with $\tau = \max\Big\{ 1+\frac{|\tau_4|}{\tau_3+\tau_4}, \frac{\tau_2|\tau_4|}{\tau_3+\tau_4}\Big\} $. Then under the same condition, we have 
\[
\|Y-E(Y)\|\leq\delta \quad \text{and}\quad \|M-E(M)\|\leq\delta
\]
holds with high probability.
\eproof

\vspace{1ex}

\noindent
{\bf Proof of Theorem  \ref{theo-spectral-init}.~~}  	
Suppose $ \tilde{\vz}_0 $ with  $ \|\tilde{\vz}_0\|=1 $ is the eigenvector of $ M $ corresponding to the largest eigenvalue $ \lambda $.  From Lemma \ref{concentration_YM}, we know
	\[
	\|M-\tau_2\|\vx\|^2I - \tau_3\vx\vx^*\|\leq \delta\|\vx\|^2.
	\]
	By Weyl's Inequality we have
	\begin{equation}\label{cond1_M}
	\big|\lambda - (\tau_2+\tau_3)\|\vx\|^2\big|\leq   \delta \|\vx\|^2.
	\end{equation}
	Then 
	\begin{align}\label{cond2_M}
 \delta\|\vx\|^2&\geq\|M-\tau_2\|\vx\|^2I-\tau_3\vx\vx^*\|\\\nonumber
	&\geq|\tilde{\vz}_0^*(M-\tau_2\|\vx\|^2I-\tau_3\vx\vx^*)\tilde{\vz}_0|\\\nonumber
	&=|\lambda-\tau_2\|\vx\|^2-\tau_3(\tilde{\vz}_0^* \vx)^2|\\\nonumber
	&\geq  \tau_3\|\vx\|^2\bigg(1-\frac{(\tilde{\vz}_0^* \vx)^2}{\|\vx\|^2}\bigg)-\big|\lambda-  (\tau_2+\tau_3)\|\vx\|^2\big|.
	\end{align}
	Substituting (\ref{cond1_M}) into (\ref{cond2_M}), we obtain
	\[
	\frac{(\tilde{\vz}_0^* \vx)^2}{\|\vx\|^2}\geq 1-\frac{2\delta}{\tau_3}.
	\]
    On the other hand, according to Lemma \ref{norm_x}, when  $ N\geq C'_\delta d $,
    \[
     \|x\|^2(1-\delta)\le\rho^2\leq(1+\delta)\,\|x\|^2
    \]
    holds with probability greater than $ 1-2\exp(-c'_\delta N) $.	
	So we can claim that for any $ \varepsilon>0 $, when $ N\geq C_{\varepsilon} d\log^2 d $ for a sufficiently large constant $ C_{\varepsilon} $, 
	\begin{align*}
	\text{dist}^2(\vz_0, \vx)&\leq \rho^2+\|\vx\|^2 - 2\rho\|\vx\|\cdot\frac{(\tilde{\vz}_0^* \vx)}{\|\vx\|}\\
	&\leq(1+\delta)\,\|\vx\|^2+\|\vx\|^2-2\sqrt{1-\delta}\,\|\vx\|^2\sqrt{1-\frac{2\delta}{\tau_3}}\\
	&\leq \bigg(3+\frac{4}{\tau_3} \bigg)\delta\,\|\vx\|^2\\
    &\leq \varepsilon^2 \|\vx\|^2
	\end{align*}
	holds with probability at least $ 1-1/d^3-2\exp(-c_{\varepsilon}\, N) $. 
\eproof

\section*{Acknowledgement}
 Yang Wang was supported in part by the Hong Kong Research Grant Council grants 16306415,  16317416 and 16308518, as well as NSFC grant 91630203. Haixia Liu is support in part by the NSFC grant 11901220. 
{\small
	\bibliographystyle{siam}
	\bibliography{phase_retrieval_fusion}
}

\section*{Appendix A. }
\setcounter{equation}{0}
\renewcommand{\theequation}{A.\arabic{equation}}
\renewcommand{\thesection}{A}

\setcounter{prop}{0}
\begin{lem}[\cite{bentkus2003inequality}]\label{one_side_concen}
 Suppose $ X_1,X_2,\ldots,X_N $ are i.i.d. real-valued random variables obeying $ X_j\leq b $ for some nonrandom $ b>0 $, $ \E X_j = 0$, and $ \E X_j^2=v^2 $. Setting $ \sigma^2= N\max(b^2, v^2) $,
 \[
 \PP(X_1+\cdots+X_N\geq y)\leq \min\Big(\exp\big(-y^2/(2\sigma^2)\big), c_0(1-\Phi\big(y/\sigma)\big)\Big)
 \]
 where one can take $ c_0 = 25 $.
\end{lem}

Now for random vectors $ \{\va_j\}_{j=1}^N $ we define the following matrix
$$
F(\vx):=\frac{1}{N}\sum^N_{j=1}\begin{bmatrix}
(A_j\vx) (A_j\vx)^* & A_j\vx( A_j\vx)^T\\
\overline{A_j\vx}( A_j\vx)^* & \overline{A_j\vx}( A_j\vx)^T
\end{bmatrix},
$$
where $A_j:=\va_j\va_j^* $. We have the following:

\begin{lem} \label{concentration_inequality3}
Let $ \{\va_j\}_{j=1}^N $ satisfying conditions (I) and (II) for generalized spectral initialization. Then 
$$
	\E\big(F(\vx)\big)
	=\begin{bmatrix}
		\tau_3\|\vx\|^2I+\tau_2\vx\vx^*+\tau_4\mathbf{D}(|\vx|^2) & (\tau_2+\tau_3)\vx\vx^T+\tau_4\mathbf{D}(\vx^2) \\
		(\tau_2+\tau_3)\overline{\vx}\vx^*+\tau_4\mathbf{D}(\overline{\vx}^2)  & \tau_3\|\vx\|^2I+\tau_2\overline{\vx}\vx^T+\tau_4\mathbf{D}(|\overline{\vx}|^2)
	\end{bmatrix},
$$
where 
\begin{align*}
\mathbf{D}(|\vx|^2) &=\diag\Big([|x_1|^2,|x_2|^2,\ldots,|x_d|^2]\Big),  \\
\mathbf{D}(\vx^2) &=\diag\Big([x_1^2,x_2^2,\ldots,x_d^2]\Big),\\
\mathbf{D}(\overline{\vx}^2) &=\diag\Big([\overline x_1^2,\overline x_2^2,\ldots,\overline{x}_d^2]\Big), \\
\mathbf{D}(|\overline{\vx}|^2) &=\diag\Big([|\overline x_1|^2,|\overline x_2|^2,\ldots,|\overline x_d|^2]\Big).
\end{align*}
For any $\delta>0$ there exist constants $ C_{\delta}, c_{\delta}>0$ such that for $ N\geq C_{\delta} \,d\log^2 d $, with probability at least $ 1-1/d^3-2\exp(-c_{\delta}\, N) $ we have
$$
\|F(\vx)-\E( F(\vx))\|\le\delta \hat{\alpha}\|\vx\|^2.
$$
Here $\hat{\alpha}:= \tau_2+\tau_3+|\tau_4|$.
\end{lem}
\Proof For any $A=[a_{ij}]$ having the same distribution as each $A_j$, by conditions (I) and (II) we have
\[
\E\big((\vx^* A\vx) A\big) = \E\Big[\Big(\sum_{i,j=1}^{d}x_j^*a_{ji}x_i\Big)\,a_{mn}\Big]=\tau_2\|\vx\|^2I+\tau_3\vx\vx^*+\tau_4\diag\Big([|x_1|^2,|x_2|^2,\ldots,|x_d|^2]\Big).
\]
It follows that $\E(a_{mm}^2) =\tau_2+ \tau_3+\tau_4$,  $\E(a_{ii}a_{mm}) = \tau_2$ for $ i\neq m $ and
$$
	\E(\overline{a_{ij}}a_{mn}) =\left\{
	\begin{array}{ll}
	\tau_3, & i=m, \, j=n\\
     0, & \text{otherwise}.
	\end{array}
	\right.
$$
Now we obtain the expectation
\[
\E\big(A\vx(A\vx)^*\big)  = \tau_3\|\vx\|^2I+\tau_2\vx\vx^*+\tau_4\mathbf{D}(|\vx|^2) 
\]
and
\[
\E\big(A\vx(A\vx)^T\big) =(\tau_2+\tau_3)\vx\vx^T+\tau_4\mathbf{D}(\vx^2). 
\]
For $ A_j = \va_j\va_j^* $,  $ \frac{1}{N}\sum_{j=1}^{N} A_j\vx(A_j\vx)^* =  \frac{1}{N}\sum_{j=1}^{N} (\vx^*A_j\vx) A_j $, we know $ \tau_2 = \tau_3 $.  And furthermore according to  Lemma  \ref{concentration_YM},         
when $ N\geq C_{\delta_1} d \log ^2d $, with probability at least $ 1-1/d^3-2\exp(-c_{\delta_1}\, N) $ we have
\[
\Big\|\frac{1}{N}\sum_{j=1}^{N}(A_j\vx)(A_j\vx)^*-\E\Big(\frac{1}{N}\sum_{j=1}^{N}(A_j\vx)(A_j\vx)^*\Big)\Big\|\leq\delta_1\hat\alpha\|\vx\|^2.
\]
While for matrix $ \frac{1}{N}\sum_{j=1}^{N} A_j\vx(A_j\vx)^T $, following the same method given in Lemma \ref{concentration_YM},  
\[
\Big\|\frac{1}{N}\sum_{j=1}^{N} A_j\vx(A_j\vx)^T -\E\Big(\frac{1}{N}\sum_{j=1}^{N} A_j\vx(A_j\vx)^T\Big)\Big\|\leq\delta_1\hat\alpha\|\vx\|^2
\]
also holds with high probability under the same condition.
The lemma is proved.
\eproof

\begin{lem} \label{lem_real}
Under the setup of Lemma \ref{concentration_inequality3}, for any $ \delta>0 $ and $\vh\in\C^d$ with $ \|\vh\|_2=1 $ satisfying $\rm Im(\vh^*\vx)=0 $, there exist constants $ C_{\delta}, c_{\delta} >0$ such that for $ N\geq C_{\delta} d \log^2 d$, with probability at least $ 1-1/d^3-2\exp(-c_{\delta}  N) $, we have
\[
	\Bigg|\frac{1}{N}\sum^N_{j=1}\Big(\re(\vh^*A_j\vx)\Big)^2-\E\bigg(\frac{1}{N}\sum^N_{j=1}\Big(\re(\vh^*A_j\vx)\Big)^2\bigg)\Bigg|\le\frac{\delta}{2}.
\]
\end{lem}
\Proof
First we have
$$
   \frac{1}{N}\sum^N_{j=1}\left(\re(\vh^*A_j\vx)\right)^2=\frac{1}{4}\begin{bmatrix}\vh\\
    \overline{\vh} \end{bmatrix}^*  F(\vx)
    \begin{bmatrix}\vh\\
     \overline{\vh}
     \end{bmatrix}.
$$
Moreover, $\|F(\vx)-\E(F(\vx))\|\le \delta\|\vx\|^2$ by Lemma \ref{concentration_inequality3}, with
\begin{eqnarray*}
	\E(F(\vx))
	&=&\begin{bmatrix}
		\tau_3\|\vx\|^2I+\tau_2\vx\vx^*+\tau_4\mathbf{D}(|\vx|^2) & (\tau_2+\tau_3)\vx\vx^T+\tau_4\mathbf{D}(\vx^2) \\
		(\tau_2+\tau_3)\overline{\vx}\vx^*+\tau_4\mathbf{D}(\overline{\vx}^2)  & \tau_3\|\vx\|^2I+\tau_2\overline{\vx}\vx^T+\tau_4\mathbf{D}(|\overline{\vx}|^2)
	\end{bmatrix}\\
	&~& \\
	&=& \tau_3\|\vx\|^2I_{2d}+\Big(\tau_2+\frac{\tau_3}{2}\Big)\begin{bmatrix}
		\vx\\ \overline{\vx}
	\end{bmatrix}\begin{bmatrix}
		\vx^* & \vx^T
	\end{bmatrix}-\frac{\tau_3}{2}\begin{bmatrix}
		\overline{\vx}\\ -\vx
	\end{bmatrix}\begin{bmatrix}
		\vx^* & \vx^T
	\end{bmatrix}+
	\tau_4\begin{bmatrix}
		\mathbf{D}(|\vx|^2) & \mathbf{D}(\vx^2)\\
		\mathbf{D}(\overline{\vx}^2)&\mathbf{D}(|\overline{\vx}|^2)
	\end{bmatrix}.
\end{eqnarray*}
Thus we have 
\begin{eqnarray*}
	&&\E\bigg(\frac{1}{N}\sum^N_{j=1}\Big(\re(\vh^*A_j\vx)\Big)^2\bigg)=\frac{1}{4}\begin{bmatrix}
		\vh\\
		\overline{\vh}
	\end{bmatrix}^*\E\big(F(\vx)\big)
	\begin{bmatrix}
		\vh\\
		\overline{\vh}
	\end{bmatrix}\\
 &&=\frac{\tau_3}{2}\|\vx\|^2\|\vh\|^2+\Big(\tau_2+\frac{\tau_3}{2}\Big)\re(\vx^*\vh)^2+\frac{\tau_4}{2}\Big(\vh^*\mathbf{D}(|\vx|^2)\vh+\re\big(\vh^*\mathbf{D}(\vx^2)\overline{\vh}\big)\Big)
\end{eqnarray*}
and  
\begin{equation*}
	\Bigg|\frac{1}{N}\sum^N_{j=1}\Big(\re(\vh^*A_j\vx)\Big)^2-\E\bigg(\frac{1}{N}\sum^N_{j=1}\Big(\re(\vh^*A_j\vx)\Big)^2\bigg)\Bigg|
     =\frac{1}{4}\begin{bmatrix}
		\vh\\
		\overline{\vh}
	\end{bmatrix}^*\Big\|F(\vx)-\E\big(F(\vx)\big)\Big\|
	\begin{bmatrix}
		\vh\\
		\overline{\vh}
	\end{bmatrix}\leq
	\frac{\delta}{2}\|\vx\|^2\|\vh\|^2.
\end{equation*}
Here the last inequality is based on Lemma \ref{concentration_inequality3}.
\eproof

\begin{lem}\label{lem_initial_bound}
	Given any fixed $ \vx $ and $ \vh $ with $ \|\vx\|_2=\|\vh\|_2=1 $ and $ \im(\vx^*\vh) =0$. Suppose $ A $ is a random matrix satisfying conditions (I) and (II) and define $ \varepsilon_0:=\frac{10}{27\alpha} \left(\sqrt{36|\tau_4|^2+\frac{27\alpha\beta}{10}} -6|\tau_4|\right)$ with $ \beta = \tau_3-(\tau_4)_- $ and $ \alpha:= \tau_2+\tau_3-(\tau_4)_-$. Then when $ s\leq \varepsilon_0 $, the following holds
	\[
	2\E(\re^2(\vh^*A\vx)) + 3s\E(\re(\vh^*A\vx)|\vh^*A\vh|)+\frac{9}{10}s^2\E(|\vh^*A\vh|^2)\geq \frac{\beta}{2}.
	\]
\end{lem}
\Proof
According to Lemma \ref{lem_real} and the condition of measurements, we easily obtain
\begin{align*}
\E(\re^2(\vh^*A\vx)) &= \frac{\tau_3}{2}+\Big(\tau_2+\frac{\tau_3}{2}\Big)\re(\vx^*\vh)^2+\frac{\tau_4}{2}\Big(\vh^*\mathbf{D}(|\vx|^2)\vh+\re\big(\vh^*\mathbf{D}(\vx^2)\overline{\vh}\big)\Big),\\
\E(\re(\vh^*A\vx)|\vh^*A\vh|) &= (\tau_2+\tau_3)\re(\vx^*\vh) + \tau_4\re(\vx^*\mathbf{D}(|\vh|^2)\vh),\\
\E(|\vh^*A\vh|^2) &= (\tau_2+\tau_3) + \tau_4\vh^*\mathbf{D}(|\vh|^2)\vh. 
\end{align*}
Note that $ \alpha = \tau_2+\tau_3-(\tau_4)_- $, $ \beta = \tau_3-(\tau_4)_-$, $ \re(\vx^*\vh)^2 = |\vx^*\vh|^2 $ and
\[
\vh^*\mathbf{D}(|\vx|^2)\vh\leq\frac{1}{2}\big( |\vx^*\vh|^2+1\big), 
\] 
we have
\begin{align*}
&2\E(\re^2(\vh^*A\vx)) + 3s\E(\re(\vh^*A\vx)|\vh^*A\vh|)+\frac{9}{10}s^2\E(|\vh^*A\vh|^2)\\
&\geq\tau_3+(2\tau_2+\tau_3)\re^2(\vx^*\vh)-(\tau_4)_-(|\vx^*\vh|^2+1)-3s(\tau_2+\tau_3+|\tau_4|)|\re(\vx^*\vh)|+\frac{9}{10}s^2\big(\tau_2+\tau_3-(\tau_4)_-\big)\\
&\geq \tau_3-(\tau_4)_-+\big(\tau_2+\tau_3-(\tau_4)_-\big)\Big(\re^2(\vx^*\vh)-3s|\re(\vx^*\vh)|+\frac{9}{4} s^2\Big)-6s|\tau_4|-\frac{27}{20}s^2\big(\tau_2+\tau_3-(\tau_4)_-\big)\\
&\geq\beta - 6s|\tau_4|-\frac{27}{20}s^2\alpha.
\end{align*}
Here we define $ \beta:= \tau_3-(\tau_4)_- $ and $ \alpha:= \tau_2+\tau_3-(\tau_4)_-$. Thus based on this inequality, when $ s\leq \frac{10}{27\alpha} \left(\sqrt{36|\tau_4|^2+\frac{27\alpha\beta}{10}} -6|\tau_4|\right) $ we have 
\[
2\E(\re^2(\vh^*A\vx)) + 3s\E(\re(\vh^*A\vx)|\vh^*A\vh|)+\frac{9}{10}s^2\E(|\vh^*A\vh|^2)\geq \frac{\beta}{2}.
\]
\eproof

\end{document}